\documentclass[12pt]{article}
\usepackage{amsmath,amssymb}
\usepackage{pstricks}
\usepackage{pst-node}
\usepackage{pst-poly}
\usepackage{multido}
\usepackage{tikz}

\newtheorem{result}{Theorem}
\newtheorem{deduce}{Corollary}

\newtheorem{define}{Definition}

\newcommand{\qed}{%
\ifmmode % if math mode, assume display: omit penalty etc.
\else \leavevmode\unskip\penalty9999 \hbox{}\nobreak\hfill \fi
\quad\hbox{\qedsymbol}}
\newcommand{\openbox}{\leavevmode \hbox to.77778em{%
\hfil\vrule
\vbox to.675em{\hrule width.6em\vfil\hrule}%
\vrule\hfil}}
\newcommand{\qedsymbol}{\openbox}

\newcommand{\showgrid}{}
\newcommand{\gridon}{\renewcommand{\showgrid}{\psset{subgriddiv=1,griddots=10,gridlabels=6pt}\psgrid}}
\gridon

\begin{document}
\begin{center}
{\LARGE  On the maximum number of edges in plane graph with fixed exterior face degree}
\end{center}

\vskip8pt

\centerline{   Niran Abbas Ali$^a$, Gek L. Chia$^{b,e}$,  \  Hazim Michman  Trao$^{c}$ \ and \   Adem Kilicman$^d$ }

\begin{center}
\itshape\small  $^{a, c, d}\/$Department of Mathematics, \\ Universiti Putra Malaysia, 43400 Serdang, Malaysia,  \\
 \vspace{1mm}
 $^{b}\/$Department of Mathematical and Actuarial Sciences, \\  Universiti Tunku Abdul Rahman, Sungai Long Campus,    Malaysia \\
\vspace{1mm}
 $^e\/$Institute of Mathematical Sciences, University of Malaya, \\ 50603 Kuala Lumpur,  Malaysia  \\
\end{center}

\begin{abstract}
A well known Euler's formula consequence's corollary in graph theory states that:
For a connected simple planar graph with $n$ vertices and $m$ edges, and girth $g$, we have  $m \leq \frac{g}{g-2}(n-2)$.
 We show that a connected simple plane graph with $n$ vertices and girth $g$, and exterior face of degree $h$ has at most $\frac{g}{g-2}(n-2)- \frac{1}{g-2}(h-g)$ edges.
 A \emph{convex hull $g$-angulation} is a connected plane graph in which the exterior face is a simple $h$-cycle and all inner faces are $g$-cycles.
 For a given  set $S$ of $n$ point in the plane having $h$ points in the boundary of its convex hull,
 we present the necessary and sufficient condition to obtain a convex hull $g$-angulation on $S$.
 We also determine the number of edges and inner faces in the convex hull $g$-angulation.
\end{abstract}

\vspace{1mm}
 \section{Introduction}

\vspace{1mm}
A \emph{topological} graph $G$ is a graph drawn in the plane, that is, its vertex set, $V (G)$, is a set of distinct points, and its edge set, $E(G)$, is a set of Jordan arcs.
A topological graph is \emph{simple} if there is no loop nor parallel edges.
The \emph{girth} of a graph $G$ is the length of a shortest cycle (if any) in $G$.
A graph is \emph{planar} if it can be embedded in the plane; a \emph{plane} graph has already been embedded in the plane.
We will refer to the regions defined by a plane graph as its \emph{faces}, the  unbounded region being called the \emph{exterior face}.
The number of edges bordering a particular face is called the \emph{degree} of the face.
To talk about the number of edges in plane graph we deal with girth $g\geq3$ and simple graphs (otherwise, $g\leq2$ and self-loop or multiple edges are available, and then there is no limit on the number of edges).

\vspace{1mm} The Euler formula for polyhedra is one of the classical results of mathematics.
Euler Polyhedron Formula stats: For any spherical polyhedron with $n$ vertices, $m$ edges, and $f$ faces, $n - m + f = 2$ (Harary's 11.1 \cite{h:refer}).

\vspace{1mm} Euler's formula has many consequence's corollaries. In particular, the following result is well known in graph theory:
For a connected simple planar graph $G$ with $n$ vertices and $m$ edges, and girth $g$, we have  $m \leq \frac{g}{g-2}(n-2)$ (Jungnickel's 1.5.3 \cite{j:refer}).
For a connected simple plane graph with $n$ vertices and girth $g$, and exterior face of degree $h\geq g$, we determine the maximum number of edges to be $\frac{g}{g-2}(n-2)-\frac{1}{g-2}(h-g)$ (Theorem \ref{1}).
That is, if $G$ contains $m$ edges where $\frac{g}{g-2}(n-2)-\frac{1}{g-2}(h-g)<m\leq \frac{g}{g-2}(n-2)$ for a given fixed $h>g$, then $G$ may be planar graph but it can not be embedded in the plane with an exterior face of degree $h$.

\vspace{1mm} Another Euler's formula consequence's corollary states:
If G is a connected simple plane graph with $n$ vertices and $m$ edges in which every face is a $g$-cycle, then $m = \frac{g}{g-2}(n-2)$ (Harary's 11.1(a) \cite{h:refer}).
Determination of the number of edges of a connected simple plane graph in which every inner face is a $g$-cycle while its exterior face is a simple cycle of degree $h\geq g$ is given by Theorem \ref{theorem 2}.

\vspace{1mm} These results are then applied to present new proofs for well known results in graph theory.

\vspace{1mm} Throughout this paper, a connected plane graph with $n$ vertices and girth $g$, and exterior face of degree $h$ will be denoted by $G_{g, h}$.
A \emph{convex hull $g$-angulation}, denoted by $H_{g,h}$, is a connected plane graph with $n$ vertices in which the exterior face is a simple cycle of degree $h$ and all inner faces are $g$-cycles.

\vspace{1mm}
 \section{Main Result}

\vspace{1mm} The number of edges in any simple plane graph depends on the number of vertices that are on its exterior face.
This is made precise in next result.

\vspace{1mm}
\begin{result}   \label{1} 
A connected simple plane graph with $n$ vertices and girth $g$, and exterior face of degree $h\geq g$ has $m\leq \frac{g}{g-2}(n-2)- \frac{1}{g-2}(h-g)$ edges.
\end{result}

\vspace{2mm}  \noindent
{\bf Proof:} Let $G_{g, h}$ be a connected simple plane graph with $n$ vertices and girth $g$, and exterior face of degree $h\geq g$.
Then $h$ is the number of boundary edges of the exterior face.
Let $m$ and $f$ be the number of edges and faces of $G_{g, h}$, respectively.
\vspace{1mm}Farther let $f_{g'}$ be the number of inner faces of degree $g'$.
Clearly $1 + \sum_{g'}{f_{g'}} =f$.

\vspace{1mm}Counting the bounding edges of all the faces we count every edge exactly twice (since each edge belongs to exactly two faces) and this follows $2m= h + \sum_{g'} g'f_{g'}$.
\vspace{1mm}Now since  $g= min\{g'\}$, that yields  $2m\geq h +  g \sum_{g'} f_{g'} = h + g(f-1)$.
Since $G_{g, h}$ is connected, we can apply Euler's formula, $n-m+f=2$.
Substituting for $f-1$ in the inequality  yields   $-2m\leq  -h - g(1-n+m)= -h -g+ng-mg $
\vspace{1mm}and then $mg-2m\leq g(n-2)- h+g$.
\vspace{1mm} Thus, $m\leq \frac{g}{g-2}(n-2)- \frac{h-g}{g-2}$.

\vspace{1mm}
\begin{define} \

(i) A convex hull $g$-angulation, $H_{g, h}$, is a connected simple plane graph with $n$ vertices in which the exterior face is a simple $h$-cycle and all inner faces are $g$-cycles.

(ii) A convex $g$-angulation, $H_{g, n}$, is a connected simple plane graph with $n$ vertices in which the exterior face is a simple $n$-cycle and all inner faces are $g$-cycles.

(iii) A $g$-angulation, $H_{g,g}$, is a connected simple plane graph with $n$ vertices in which all faces are $g$-cycles.
\end{define}

\vspace{1mm}
\begin{result}  \label{theorem 2}
Let $S$ be a set of $n$ points in the plane having $h$ vertices in the convex hull and let $g\geq3$.  A simple plane graph on $S$ with girth $g$ and exterior face of degree $h\geq g$, is a connected convex hull $g$-angulation $H_{g, h}$  if and only if  $2n-h-g$ is divisible by  $g-2$ and $m=\frac{g}{g-2}(n-2)- \frac{1}{g-2}(h-g)$.
Moreover, $H_{g, h}$ has exactly  $m=n+t$ edges and $t+1$ inner faces where $t=\frac{2n-h-g}{g-2}$.
\end{result}

\vspace{1mm}  \noindent
{\bf Proof:} Let $S$ be a set of $n$ points in the plane having $h$ vertices in the boundary of its convex hull and let $g\geq3$.

\vspace{1mm} Assume that $G_{g, h}$ be a connected simple plane graph on $S$ with  $m=\frac{g}{g-2}(n-2)- \frac{1}{g-2}(h-g)$ edges and girth $g$, and exterior face of degree $h\geq g$. Assume further that $2n-h-g$ is divisible by  $g-2$.

\vspace{1mm} Suppose on the contrary that there is an inner face of degree $g'>g$. Assume without loss of generality that $g'=g+1$.
By repeating the same argument of previous proof, we see that $2m= g(f-2)+ h +g+1$.
Since $G_{g, h}$ is connected, we can apply Euler's formula, $n-m+f=2$.
Substituting for $f-2$  yields   $2m= g(m-n)+ h +g+1$,      and then $mg-2m=m g(n-2)- h+g -1$     which follows $m=\frac{g}{g-2}(n-2)- \frac{1}{g-2}(h-g)-\frac{1}{g-2}$, a contradiction. Hence, $G_{g, h}$ is a $g$-angulation.

\vspace{2mm}  Assume that, $H_{g, h}$ is a convex hull $g$-angulation.
Since each inner face of $H_{g, h}$  is of degree $g$, then by repeating the same argument of previous proof, we see that $2m=  h + g f_{g} = h + g (f-1)$.
Substituting for $f-1$ of Euler’s Formula ($n-m+f=2$) yields $ m=\frac{g}{g-2}(n-2)- \frac{1}{g-2}(h-g)$.
Now, since $m$ is a natural number then  $g(n-2)-(h-g)=2n-h-g+n(g-2)$ is divisible by $g-2$.
Hence, $2n-h-g$ is divisible by $g-2$.
Let $t=\frac{2n-h-g}{g-2}$. Thus, $m=n+t$.

\vspace{1mm} Now, by (Euler Polyhedron Formula), the number of inner faces is $f-1=m-n+1=t+1$. \qed

\

\vspace{2mm} As a direct consequence of Theorem \ref{1} and Theorem \ref{theorem 2}, we have the following Corollaries.

\vspace{1mm}
\begin{deduce}     \label{corollary 4}
A convex connected simple plane graph  $G_{g, n}$ with $n$ vertices and girth $g$ has  $m\leq \frac{(g-1)n  - g}{g-2}$ edges.
\end{deduce}

\vspace{1mm}  \noindent
{\bf Proof:}  By applying Theorem \ref{1} with inserting $h=n$.  \qed

\vspace{1mm}
\begin{deduce}     \label{corollary 5}
A convex $g$-angulation $H_{g,n}$ with $n$ vertices where $n=g+t(g-2)$ for some integer $t\geq0$  has $m = \frac{(g-1)n  - g}{g-2}$ edges. Moreover, $m=n+t$ and the number of inner faces is $t+1$.
\end{deduce}

\vspace{1mm}  \noindent
{\bf Proof:}  Assume that $H_{g,n}$ is a convex $g$-angulation with $n=g+t(g-2)$ vertices for some integer $t\geq0$.
By applying Theorem \ref{theorem 2} with substituting for $h=n$ yields $m=\frac{(g-1)n  - g}{g-2}$ and  $\frac{2n-h-g}{g-2}=\frac{n-g}{g-2}=t$.
 Hence, $m=n+t$ and $f-1=t+1$.  \qed

\

\vspace{2mm}

An important result in planarity that any outerplane graph (or convex connected simple plane graph) $G_{3, n}$ of girth $3$ has $ m\leq 2n-3$ edges and the equality is achieved when the graph is a maximal outerplane (or a convex triangulation) $H_{3,n}$ can be proved by applying Corollary \ref{corollary 4} and Corollary \ref{corollary 5} with substituting for $g=3$.

For following well known results in graph theory, we present new proofs depending on Theorem \ref{1} and Theorem \ref{theorem 2}:

\vspace{1mm}
\begin{result}    \label{corollary 1}
(9.1 \cite{bcko:refer})
Any convex hull triangulation, of $n$ vertices having exterior face of degree $h$, has $3n - 3 - h$ edges and $2n-2-h$ inner triangles.
\end{result}

\vspace{1mm}  \noindent
{\bf Proof:}  By applying Theorem \ref{theorem 2} with substituting for $g=3$ yields $t=2n-3-h$.
Hence, $m =n+t =3m-3-h $ and the number of inner faces is $t+1 =2n-2-h$. \qed

\vspace{1mm}
\begin{result}    \label{corollary 2}
(Jungnickel's 1.5.3 \cite{j:refer})
A connected simple plane graph $G_{g, g}$ with $n$ vertices and girth $g$, and exterior face of degree $g$ has $m\leq \frac{g}{g-2}(n-2)$ edges.
\end{result}

\vspace{1mm}  \noindent
{\bf Proof:}  By applying Theorem \ref{1}  
with substituting for $h=g$.  \qed

\vspace{1mm}
\begin{deduce}     \label{corollary 3}
(Harary's 11.1(a) \cite{h:refer})
A $g$-angulation $G_{g,g}$ with $n$ vertices, where $n=g+t'(g-2)$ for some integer $t'\geq0$, has $m= \frac{g}{g-2}(n-2)$ edges. Moreover, $m=n+2t'$ and the number of inner faces is $2t'+1$.
\end{deduce}

\vspace{1mm}  \noindent
{\bf Proof:}  Assume that $G_{g,g}$ is a $g$-angulation with $n=g+t'(g-2)$  vertices for some integer $t'\geq0$.
By applying Theorem \ref{theorem 2} with substituting for $h=g$ yields $m = \frac{g}{g-2}(n-2)$ and  $t=\frac{2n-h-g}{g-2}=2(\frac{n-g}{g-2})=2t'$.
Hence, $m=n+2t'$ and the number of inner faces is $f-1=2t'+1$.  \qed

\end{document}